\documentclass[12pt]{amsart}
\usepackage{latexsym,amsmath,amscd}

\newtheorem{theorem}{Theorem}

\newtheorem{proposition}{Proposition}

\newtheorem{lemma}{Lemma}
\theoremstyle{remark}

\title{an arakelov theoretic proof of the equality of conductor and discriminant}
\date{May 30, 2000}
\author{s\.{i}nan \"{u}nver}
\address{Department of Mathematics\\ University of California-Berkeley\\Berkeley, CA 94720}
\email{sinan@math.berkeley.edu}
\begin{document}
\maketitle
\noindent

\section{Introduction}
Let $K$ be a number field, $\mathcal{O}_K$ be the ring of integers of $K$, and 
$S$ be $\hbox{Spec}(\mathcal{O}_K)$. Let $f:X \to S$ be an arithmetic surface. By this we mean a regular scheme, proper 
and flat over $S$, of relative 
dimension one. We also assume that the generic fiber of $X$ has genus$\ge1$, and that $X/S$ has geometrically connected fibers.     
 
 Let $\omega_{X}$ be the dualizing sheaf of $X/S$. 
The Mumford isomorphism ([Mumf], Theorem 5.10) \[\det \hbox{R}f_{*}(\omega_{X}^{\otimes2})\otimes K \to (\det \hbox{R} f_{*}\omega_{X})^{\otimes13}\otimes K,\] which is unique up to sign, 
gives a rational section $\Delta$ of \[(\det \hbox{R}f_{*}\omega_{X})^{\otimes13}\otimes(\det \hbox{R} f_{*}(\omega_{X}^{\otimes2}))^{\otimes-1}.\] The discriminant $\Delta(X)$ 
of $X/S$ is defined as the divisor of this rational section ([Saito]). If $\mathfrak{p}$ is a closed point of $S$, we denote the coefficient 
of $\mathfrak{p}$ in $\Delta(X)$ by $\delta_{\mathfrak{p}}$. 

On the other hand $X/S$ has an Artin conductor $\hbox{Art}(X)$ (cf. [Bloch]), which is similarly a divisor on
$S$. We denote the coefficient of $\mathfrak{p}$ in $\hbox{Art}(X)$ by $\hbox{Art}_{\mathfrak{p}}$. 
Let $S'$ be the strict henselization of complete local ring 
at $\mathfrak{p}$, with field of fractions $K'$. Let s be its special point, $\eta$ be its generic point, and $\overline{\eta}$ be a
geometric  generic point corresponding to an algebraic closure $\overline{K'}$ of $K'$. Let $\ell$ be a prime different from the residue characteristic at $\mathfrak{p}$. 
Then 
\begin{eqnarray*} 
\hbox{Art}_{\mathfrak{p}}(X)&=&\sum_{i \geq 0} {(-1)^{i}
\dim_{\mathbb{Q}_{\ell}}\hbox{H}^{i}_{\acute{e}t}
(X_{\overline{\eta}},\mathbb{Q}_{\ell})}
-\sum_{i \geq 0}{(-1)^{i} \dim_{\mathbb{Q}_{\ell}}
\hbox{H}^{i}_{\acute{e}t}(X_{s},\mathbb{Q}_{\ell})}\\
&&+\sum_{i \geq 0}{(-1)^{i} \hbox{Sw}_{\overline{K'}/K'}
(\hbox{H}^{i}_{\acute{e}t}(X_{\overline{\eta}},\mathbb{Q}_{\ell}))},
\end{eqnarray*}
where $\hbox{Sw}_{\overline{K'}/K'}$ denotes the Swan conductor of the Galois representation of $\overline{K'}/K'$.
Both of these divisors are supported on the primes of bad reduction of $X$. 
We give another proof of Saito's theorem ([Saito], Theorem 1) in the number field  case.

\begin{theorem}
For any closed point $\mathfrak{p}\in S$, we have $\delta_{\mathfrak{p}}=-\hbox{Art}_{\mathfrak{p}}.$
\end{theorem}
Fix a K\"{a}hler metric on $X$, this gives metrics on $\Omega^{1}_{X_{\nu}}$'s, 
for each $\nu\in{S(\mathbb{C})}$. For a hermitian coherent sheaf $\mathcal{E}$, we endow 
$\det\hbox{R}f_{*}(\mathcal{E})$ with its Quillen metric. The proof of the theorem has the following corollaries. 

\begin{proposition}
We have 
\begin{eqnarray*}
\deg \det \hbox{R}f_{*}\omega_{X}
&=& \frac{1}{12}[\deg f_{*}(\widehat{c_{1}}(\omega_{X})^{2})
+\log Norm(-\hbox{Art}(X))]\\
&&  [K:\mathbb{Q}](g-1)[2\zeta '(-1)+\zeta(-1)],   
\end{eqnarray*}
with $\zeta$ the Riemann zeta function.
\end{proposition} 
Proposition 1 is an arithmetic analogue of Noether's formula in which $\det 
\hbox{R}f_{*}\omega_{X}$ is endowed with the Quillen metric. 
Faltings [Falt] and Moret-Bailly [M-B] proved 
a similar formula for the Faltings metrics. 

\begin{proposition}
We have 
\[\frac{1}{[K:\mathbb{Q}]}\sum_{\nu \in S(\mathbb{C})}{\log \Vert \Delta_{\nu} \Vert}=12(1-g)[2\zeta '(-1)+\zeta(-1)].\]
In particular, the norm of the  Mumford isomorphism does not depend on the metric. 
\end{proposition}
\section{Proof}

First we prove Proposition 1. 
By duality ([Deligne], Lemme 1.3), $\deg \det \hbox{R}f_{*}\omega_{X}=\deg \det \hbox{R}f_{*}\mathcal{O}_{X}$. 
 By the arithmetic Riemann-Roch theorem of Gillet and Soul\'{e} ([G-S], Theorem 7),we get  
\[ \deg \det \hbox{R}f_{*}\mathcal{O}_{X}=\deg f_{*}(\widehat{Td}(\Omega^{1}_{X})^{(2)})-\frac{1}{2}\sum_{\nu\in{S(\mathbb{C})}}\int_{X_{\nu}}Td(T_{X_{\nu}})R(T_{X_{\nu}})\]
Here $Td$ and $R$ are the Todd and Gillet-Soul\'{e} genera respectively, and the upperscript
$(2)$ denotes the degree $2$ component. Applying the definitions of these characteristic classes
we obtain \[\deg \det \hbox{R}f_{*}\mathcal{O}_{X}=\frac{1}{12}\deg f_{*}(\widehat{c_{1}}(\Omega_{X}^{1})^{2}+\widehat{c_{2}}(\Omega_{X}^{1}))+[K:\mathbb{Q}](g-1)[2\zeta'(-1)+\zeta(-1)]\]
Let $Z$ denote the union of singular fibers of $f$, and  let $c_{2,X}^{Z}(\Omega^{1}_{X})$ be the localized Chern class of $\Omega^{1}_{X}$ with support in $Z$ (cf. [Bloch], [Fulton]).  
 Chinburg, Pappas, and Taylor ([CPT], Proposition 3.1) prove the formula \[\deg f_{*}(\widehat{c_{2}}
(\Omega^{1}_{X}))=\log Norm(c_{2,X}^{Z}(\Omega^{1}_{X})).\] 
 Combining this with the fundamental formula of Bloch ([Bloch], Theorem 1)  
 \[-\hbox{Art}_{\mathfrak{p}}(X)=\deg_{\mathfrak{p}} c_{2,X}^{Z}(\Omega_{X}^{1}),\]  we obtain the desired 
formula. Note that, since $\det\Omega_{X}^{1}=\omega_{X}$,  $\widehat{c_{1}}(\Omega_{X}^{1})=\widehat{c_{1}}(\omega_{X})$.
\hfill $\Box$
Taking degrees in the Mumford isomorphism gives \[13\deg\det\hbox{R}f_{*}\omega_{X}
=\deg\det\hbox{R}f_*(\omega_{X}^{\otimes2})+\log Norm(\Delta(X))-\sum_{\nu \in S(\mathbb{C})}\log \Vert \Delta_{\nu} \Vert.\]
The arithmetic Riemann-Roch theorem gives 
 
\[\deg\det\hbox{R}f_{*}(\omega_{X}^{\otimes2})=\deg\det\hbox{R}f_{*}\omega_{X}+\deg f_{*}(\widehat{c_{1}}(\omega_{X})^{2}).\]
Therefore we get 
\begin{equation}
\deg\det\hbox{R}f_{*}\omega_{X}
=\frac{1}{12}[\deg f_{*}(\widehat{c_{1}}(\omega_{X})^{2})
+\log Norm(\Delta(X))-\sum_{\nu \in S(\mathbb{C})}{\log \Vert \Delta_{\nu} \Vert }.
\end{equation}
Subtracting (1) from the expression in the statement of Proposition 1, we obtain
\begin{equation}
\log (\frac{Norm(\Delta(X/S))}{Norm(-\hbox{Art}(X/S))})= \sum_{\nu \in S(\mathbb{C})}{\log \Vert \Delta_{\nu} \Vert}+12[K:\mathbb{Q}](g-1)[2\zeta'(-1)+\zeta(-1)].
\end{equation}

Now $X_{K}$ has semistable reduction after a finite base change $K'/K$. For semistable 
 $X'/S'$, both $-\hbox{Art}_{\mathfrak{p}'}(X')$ ([Bloch]), and $\delta_{\mathfrak{p}'}(X')$ 
([Falt], Theorem 6) are equal to the number of singular points in the geometric fiber over $\mathfrak{p}'$. 
Therefore $-\hbox{Art}_{\mathfrak{p}'}=\delta_{\mathfrak{p}'}$, and hence 
\begin{equation}
Norm(\Delta(X'/S'))=Norm(-\hbox{Art}(X'/S')).
\end{equation}
 Applying this 
to a semistable model $X'$ of $X\otimes_{K} {K'}$, 
and noting that the base change multiplies the right hand side of (2) by $[K':K]$, we see that 
the right hand side of (2) is equal to zero, and hence that the equality  
 
\begin{equation}
Norm(\Delta(X/S))=Norm(-\hbox{Art}(X/S)) 
\end{equation}
holds for $X$.

To prove the equality $\delta_{\mathfrak{p}}=-\hbox{Art}_{\mathfrak{p}}$ for an arbitrary 
closed point $\mathfrak{p}\in S$, we will use the following lemma.
\begin{lemma}
Fix distinct closed points $\beta_{1},..,\beta_{s}\in S$. And fix finite extensions 
$L_{i}$ of the completions $K_{i}$ of $K$ at $\beta_{i}$'s, for each $1 \leq i \leq s$,  such that 
$[L_{i}:K_{i}]=n$ for some $n$. Then there exists an extension $L/K$ such that, for each $1 \leq i \leq s$, there
is only one prime $\gamma_{i}$ of $L$ lying over $\beta_{i}$, and the completion of $L$ at $\gamma_{i}$ is isomorphic (over $K_{i}$) to 
$L_{i}$.
\end{lemma} 

Proof. The proof is an easy application of Krasner's lemma, and the approximation lemma. 
Details are omitted.
\hfill $\Box$

Take $\mathfrak{p}$=$\beta_{1}$, a prime of bad reduction. Denote the remaining  primes 
of bad reduction by $\beta_{i}$, $2 \leq i \leq s$. Choose  extensions $L_{i}$ of the 
local fields $K_{i}$, for all $1 \leq i \leq s$, such that $L_{1}$ is unramified over $K_{1}$, $X$ has semistable reduction 
over $L_{i}$, for $2 \leq i \leq s$, and $[L_{i}:K_{i}]=n$, for some n. Applying the lemma to this data 
we obtain an extension $L$ of $K$. Let $T=\hbox{Spec}(\mathcal{O}_{L})$. The curve $X\otimes_{K}L$ has a proper, regular model $Y$ over $T$ such that
\begin{itemize}
\item[(i)] $ Y\otimes_{T} T_{\gamma_{1}} \simeq X\otimes_{S} T_{\gamma_{1}}$, and 
\item[(ii)]  $Y$ is semistable at $\gamma_{i}$, for $2 \leq i \leq s$.
\end{itemize}   
Applying (4) to $Y$ gives the equality
\[\sum_{1 \leq i \leq s}{\delta_{\gamma_{i}}\log Norm(\gamma_{i})}=\sum_{1 \leq i \leq s}{ -\hbox{Art}_{\gamma_{i}}\log Norm(\gamma_{i})}.\] On the other hand  
because of semistability, we have  $\delta_{\gamma_{i}}=-\hbox{Art}_{\gamma_{i}}$, for $ 2 \leq i \leq s$. Hence we get $\delta_{\gamma_{1}}=-\hbox{Art}_{\gamma_{1}}$. Since 
 $T/S$ is \'{e}tale at $\gamma_{1}$ , (i) implies \[\delta_{\mathfrak{p}}=\delta_{\gamma_{1}}=-\hbox{Art}_{\gamma_{1}}=-\hbox{Art}_{\mathfrak{p}}.\]    
 
 Acknowledgements. I would like to thank A.Abbes for his many mathematical suggestions, and 
 for his constant encouragement. I would like to thank my adviser 
 P.Vojta for supporting me this academic year.

\end{document}